\documentclass[a4paper, 12pt]{article}
\usepackage{amssymb}

\usepackage{amsmath,amsthm}

\begin{document}

\title{Coordination of Multiple Dynamic Agents with Asymmetric Interactions \thanks{This work was supported by the National Natural Science Foundation of China (No. 10372002 and No.
60274001) and the National Key Basic Research and Development
Program (No. 2002CB312200).}}
\author{Hong Shi, Long Wang, Tianguang Chu\\
Intelligent Control Laboratory, Center for Systems and Control,\\
Department of Mechanics and Engineering Science,\\
Peking University, Beijing 100871, P. R. China}
\date{}
\maketitle


\textbf{Abstract}-- In this paper, we consider multiple mobile
agents moving in Euclidean space with point mass dynamics. Using a
coordination control scheme, we can make the group generate stable
flocking motion. The control laws are a combination of
attractive/repulsive and alignment forces, and the control law
acting on each agent relies on the position information of all
agents in the group and the velocity information of its neighbors.
By using the control laws, all agent velocities become
asymptotically the same, collisions can be avoided between all
agents, and the final tight formation minimizes all agent global
potentials. Moreover, we show that the velocity of the center of
mass is invariant and is equal to the final common velocity.
Furthermore, we study the motion of the group when the velocity
damping is taken into account. We prove that the common velocity
asymptotically approaches zero, and the final configuration
minimizes the global potential of all agents. In this case, we can
properly modify the control scheme to generate the same stable
flocking. Finally, we provide some numerical simulations to
further illustrate our results.

\textbf{Keywords}---Collective behavior, swarms, robot teams,
coordination, flocking, asymmetric interactions, multi-agent
systems, collision avoidance, stability.


\section{Introduction}

In nature, flocking can be found everywhere and it can be regarded
as a typical behavior of large number of interacting dynamic
agents. This exists in the form of flocking of birds, schooling of
fish, and swarming of bacteria. Understanding the mechanisms and
operational principles in them can provide useful ideas for
developing distributed cooperative control and coordination of
multiple mobile autonomous agents/robots. In recent years,
distributed control/coordination of the motion of multiple dynamic
agents/robots has emerged as a topic of major interest \cite{N. E.
Leonard and E. Fiorelli}--\cite{John H. Reif; H. Wang}. This is
partly due to recent technological advances in communication and
computation, and wide applications of multi-agent systems in many
engineering areas including cooperative control of unmanned aerial
vehicles (UAVs), scheduling of automated highway systems,
schooling for underwater vehicles, attitude alignment for
satellite clusters and congestion control in communication
networks \cite{F. Giulietti}--\cite{R. Olfati-Saber 1}.
Correspondingly, there has been considerable effort in modelling
and exploring the collective dynamics in physics, biology, and
control engineering, and trying to understand how a group of
autonomous creatures or man-made mobile autonomous agents/robots
can cluster in formations without centralized coordination and
control \cite{T. Vicsek}--\cite{T. Chu3}.

In order to generate computer animation of the motion of flocks,
Reynolds \cite{C. W. Reynolds} modelled the boid as an object
moving in a three dimensional environment based on the positions
and velocities of its nearby flockmates and introduced the
following three rules (named steering forces) \cite{C. W.
Reynolds}:

1) Collision Avoidance: avoid collisions with nearby flockmates,

2) Velocity Matching: attempt to match velocity with nearby
flockmates, and

3) Flock Centering: attempt to stay close to nearby flcokmates.\\
Subsequently, Vicsek {\it et al.} \cite{T. Vicsek} proposed a
simple model of autonomous agents (i.e., points or particles). In
the model, all agents move at a constant identical speed and each
agent updates its heading as the average of the heading of agent
itself with its nearest neighbors plus some additive noise. They
demonstrated numerically that all agents will eventually move in
the same direction, despite the absence of centralized
coordination and control. In fact, Vicsek's model can be seen as a
special case of Reynolds's model, and it only considers the
velocity matching between agents. Jadbabie {\it et al.} \cite{A.
Jadbabaie} and Savkin \cite{A. V. Savkin} used two kinds of
completely different methods to provide the theoretical
explanation for the observed behaviors in Vicsek's model,
respectively. According to the results in \cite{C. W. Reynolds},
Tanner {\it et al} \cite{tanner 1} studied a swarm model that
consists of multiple mobile agents moving on the plane with double
integrator dynamics. They introduced a set of control laws that
enabled the group to generate stable flocking motion and provided
strictly theoretical justification. However, it is perhaps more
reasonable to take the agents' masses into account and consider
the point mass model in which each agent moves in $n$-dimensional
space based on the Newton's law. In this paper, we investigate the
collective behavior of multi-agent systems in $n$-dimensional
space with point mass dynamics.

In \cite{tanner 1}, the authors used an undirected graph to
describe the neighboring relations between agents, which means
that the neighboring relations are mutual. In other words, they
only considered the case with bidirectional information exchange
between agents. However, under some circumstances, the information
exchange is not mutual. In fact, due to the agent differences,
they maybe have different action forces on different agents and
even have different sense ranges, hence, the influence intensities
between two agents might be different with each other and even
their information can not be exchanged with each other at all. For
example, in a group of agents with spherical sense neighborhoods
but with different radii of the neighborhoods or a group of agents
with conic sense neighborhoods, the information exchange among
them might be unidirectional. A group of mobile robots with conic
vision range is just an example. In this paper, the results in
\cite{tanner 1} are extended to a directed graph. We consider the
stability properties of the group in the case of directed
information exchange. In order to generate stable flocking, we
introduce a set of control laws so that each agent regulates its
velocity based on a fixed set of ``neighbors" and regulates its
position such that its global potential become minimum. Note that,
in this paper, we only consider the fixed topology of the
neighboring relations, and the case that the information topology
is dynamic will be discussed in another paper. Here, the control
laws are a combination of attractive/repulsive and alignment
forces. By using the control laws, all agent velocities become
asymptotically the same, collisions can be avoided between all
agents, and the final tight formation minimizes all agent
potentials.

This paper is organized as follows: In Section 2, we formulate the
problem to be investigated. Some basic concepts and results in
graph theory are provided in Section 3. We analyze the system
stability with some specific control laws in Section 4. Some
numerical simulations are presented to further illustrate our
results in Section 5. Finally, we briefly summarize our results in
Section 6.

\section{Problem Formulation}

We consider a group of $N$ agents moving in an $n$-dimensional
Euclidean space, each has point mass dynamics described by
\begin{equation}
\begin{array}{l}
\quad\,\dot{x}^i=v^i, \\
m_i\dot{v}^i=u^i,\ \ i=1, \cdots, N, \\
\end{array}
\label{eq1}
\end{equation}
where $x^i=(x_{1}^i, \cdots, x_{n}^i)^T\in R^n$ is the position
vector of agent $i$; $v^i=(v_{1}^i, \cdots, v_{n}^i)^T\in R^n$ is
its velocity vector, $m_i>0$ is its mass, and $u^i=(u_{1}^i,
\cdots, u_{n}^i)^T\in R^n$ is the (force) control input acting on
agent $i$. $x^{ij}=x^i-x^j$ denotes the relative position vector
between agents $i$ and $j$.

Our aim is to make the whole group move at a common velocity and
maintain constant distances between all agents. We first consider
the ideal case, that is, we ignore the velocity damping. In order
to achieve our objective, we try to decrease the velocity
differences between agents, and at the same time, regulate their
distances such that their global potentials become minimum. Hence,
we choose the control law for each agent to be a combination of
two components. The control input $u^i$ for agent $i$ is
\begin{equation}
u^i=\alpha^i+\beta^i, \label{eq2}
\end{equation}
where $\alpha^i$ is used to regulate the potentials among agents
and $\beta^i$ is used to regulate the velocity of agent $i$ to the
weighted average of its ``neighbors". $\alpha^i$ is derived from
the social potential fields which is described by artificial
social potential function, $V^{i}$, which is a function of the
relative distances between agent $i$ and its flockmates.
Collision-free and cohesion in the group can be guaranteed by this
term. Note that $\alpha^i$ indicates the tendency of collision
avoidance and cohesion of the flocks, whereas $\beta^i$ indicates
the tendency of agent velocity matching.

Certainly, in some cases, the velocity damping can not be ignored.
For example, the objects moving in viscous environment and the
mobile objects with high speeds, such as air vehicles, are subject
to the influence of velocity damping. Then, under these
circumstances, the model in (\ref{eq1}) should be the following
form
\begin{equation}
\begin{array}{lcl}
\quad\,\dot{x}^i=v^i, \\
m_i\dot{v}^i=u^i-k_iv^i,\\
\end{array}
\label{eq3}
\end{equation}
where $k_i>0$ is the ``velocity damping gain", $-k_iv^i$ is the
velocity damping term, and $u^i$ is the control input for agent
$i$. Note that we assume the damping force is in proportion to the
magnitude of velocity. And, because the ``velocity damping gain"
is determined by the shape and size of the object, the property of
medium, and some other factors, we assume that the damping gains
$k_i$, $i=1, \cdots, N$ are not equal to each other. Certainly, in
some cases, the assumption of the same gain is enough. In order to
achieve our aim, the velocity damping should be cancelled by some
terms in the control laws. Thus, we modify the control scheme to
be
\begin{equation}
u^i=\alpha^i+\beta^i+k_iv^i. \label{eq4}
\end{equation}

\section{Main Results}

In this section, we investigate the stability properties of
multiple mobile agents with point mass dynamics described in
(\ref{eq1}). We present explicit control input in (\ref{eq2}) for
the terms $\alpha^i$ and $\beta^i$. In this paper, the control law
acting on each agent is based on two kinds of information
topologies that is the position information topology and the
velocity information topology. We will employ algebraic graph
theory as basic tools to study the properties of the group. Some
concepts and results in graph theory are given in the Appendix.

In this paper, we assume that each agent is equipped with two
onboard sensors: the position sensor which is used to sense the
position information of the flockmates and the velocity sensor
which is used to sense the velocity information of its neighbors,
and assume that all the sensors can sense instantaneously.
Correspondingly, we define two kinds of structure topologies to
describe the neighboring relations between the agents. We will use
an undirected graph ${\cal{G}}$ to describe the position sensor
information flow and use a weighted directed graph ${\cal{D}}$ to
describe the velocity sensor information flow.

First, we make the following definitions and assumptions.

\textbf{Definition 1}: (Position neighboring graph) The position
neighboring graph, $\cal{G}=(\cal{V}, \cal{E})$, is an undirected
graph consisting of a set of vertices, ${\cal{V}}=\{n_1, \cdots,
n_N\}$, indexed by the agents in the group, and a set of edges,
${\cal{E}}=\{(n_i, n_j)\in{\cal{V}\times \cal{V}}\ |\ n_j\sim
n_i\}$, which contain unordered pairs of vertices that represent
the position neighboring relations.

\textbf{Definition 2}: (Velocity neighboring graph) The velocity
neighboring graph, $\cal{D}=(\cal{V}, \cal{E})$, is a directed
graph consisting of a set of vertices, ${\cal{V}}=\{n_1, \cdots,
n_N\}$, indexed by the agents in the group, and a set of arcs,
${\cal{E}}=\{(n_i, n_j)\in{\cal{V}\times \cal{V}}\ |\ n_j\sim
n_i\}$, which contain ordered pairs of vertices that represent the
velocity neighboring relations.

Note that, in ${\cal{E}}$, an arc $(n_i, n_j)$ represents a
unidirectional velocity information exchange link from $n_i$ to
$n_j$, which means that agent $i$ can sense the velocity of agent
$j$.

\textbf{Assumption 1}: The position neighboring graph $\cal{G}$ is
complete.

In order to make the final potential of each agent be global
minimum and at the same time, ensure collision-free in the group,
we assume that the position neighboring graph is complete. This
means that, each agent can always obtain the position information
of all the other agents in the group. Certainly, in the case that
the position neighboring relation is determined by a certain
neighborhood around the agent and consequently cause the topology
of the neighboring graph ${\cal{G}}$ to be dynamic, we can also
guarantee collision avoidance in the group.

\textbf{Assumption 2}: The velocity neighboring graph $\cal{D}$ is
weakly connected.

In this paper, we consider a group of mobile agents with fixed
topology, so $\cal{D}$ is weakly connected and does not change
with time. Denote the set ${\cal{N}}_i\triangleq \{j\ |\ a_{ij}>0
\}\subseteq \{1, \cdots, N\}\backslash \{i\}$ which contains all
neighbors of agent $i$. If agent $j$ is a neighbor of agent $i$,
we denote $j\sim i$, and otherwise we denote $j\nsim i$.

\textbf{Definition 3} \cite{tanner 1}: (Potential function)
Potential $V^{ij}$ is a differentiable, nonnegative, radially
unbounded function of the distance $\|x^{ij}\|$ between agents $i$
and $j$, such that

$i$) $V^{ij}(\|x^{ij}\|)\rightarrow \infty$ as
$\|x^{ij}\|\rightarrow 0$,

$ii$) $V^{ij}$ attains its unique minimum when agents $i$ and $j$
are located at a desired distance.

Functions $V^{ij}$, $i, j=1, \cdots, N$ are the artificial social
potential functions that govern the interindividual interactions.
Cohesion and separation can be achieved by artificial potential
fields \cite{E. Rimon}. One example of such potential function is
the following
\[V^{*}(x)=a\ln{x^2}+\frac{b}{x^2},\]
where $x\in R_+=(0,\infty)$ is variable, $a>0$ and $b>0$ are some
constants. It is easy to see that $V^*$ attains its unique minimum
when $x=\sqrt{b/a}$. Hence, when the distance $\|x^i-x^j\|$
between agents $i$ and $j$ is $\sqrt{b/a}$, the potential function
$V^{ij}$ attains its unique minimum.

By the definition of $V^{ij}$, the total potential of agent $i$
can be expressed as
\begin{equation}
V^i=\sum_{j=1, j\neq i}^NV^{ij}(\|x^{ij}\|).
 \label{eq5}
\end{equation}

Agent dynamics are different in ideal case (i.e., velocity damping
is ignored) and nonideal case. This means that the agent has
different motion equations in the two cases. Hence, in what
follows, we will discuss the motion of the group in the two
different cases, respectively.

\subsection{Ideal Case}

In this case, in order to achieve our control aim, we take the
control law $u^i$ to be
\begin{equation}
u^i=-\sum_{j\in {\cal{N}}_i}w_{ij}(v^i-v^j)-\sum_{j=1, j\neq
i}^N\nabla_{x^i}V^{ij}.
 \label{eq6}
\end{equation}
Note that, $w_{ij}\geq 0$, and $w_{ii}=0$, $i, j=1, \cdots, N$
represent the interaction coefficients. And $w_{ij}>0$ if agent
$j$ is a neighbor of agent $i$, and is 0 otherwise. We denote
$W=[w_{ij}]$. Thus, by the weakly connectivity of the velocity
neighboring graph, $W+W^T$ is irreducible. The control law in
(\ref{eq6}) implies that we adopt the local velocity regulation
and the global potential regulation to achieve our aim.

In the discussion to follow, we will need the concept of weight
balance condition defined below:

\noindent \textbf{Weight Balance Condition} \cite{Wang}: consider
the weight matrix $W=[w_{ij}]\in R^{N\times N}$, for all
$i=1,\cdots,N$, we assume that
$\sum_{j=1}^Nw_{ij}=\sum_{j=1}^Nw_{ji}$.

The weight balance condition has a graphical interpretation:
consider the directed graph associated with a matrix, weight
balance means that, for any node in this graph, the weight sum of
all incoming edges equals the weight sum of all outgoing edges
\cite{R. Horn and C. R. Johnson}. The weight balance condition can
find physical interpretations in engineering systems such as water
flow, electrical current, and traffic systems.

\textbf{Proposition 1}: Let ${\cal{D}}$ be a weighted directed
graph such that the weight balance condition is satisfied. Then
${\cal{D}}$ is strongly connected if and only if it is weakly
connected.

\textbf{Proof}: It is obvious that if ${\cal{D}}$ is strongly
connected, then it is weakly connected. Hence, we only need to
prove that if ${\cal{D}}$ is weakly connected, then it is strongly
connected. In the following, we will use the way of contradiction
to prove it. Assume that ${\cal{D}}$ is weakly connected, but not
strongly connected, then we denote all strongly connected
components of ${\cal{D}}$ as ${\cal{D}}_1, \cdots, {\cal{D}}_m$,
where $m$ is an integer and $m>1$. If there is an arc starting in
${\cal{D}}_i$ and ending in ${\cal{D}}_j$, then any arc joining
${\cal{D}}_i$ to ${\cal{D}}_j$ must start in ${\cal{D}}_i$. Hence
we can define a directed graph ${\cal{D}}^*$ with the strongly
connected components of ${\cal{D}}$ as its vertices, and such that
there is an arc from ${\cal{D}}_i$ to ${\cal{D}}_j$ in
${\cal{D}}^*$ if and only if there is an arc in ${\cal{D}}$
starting in ${\cal{D}}_i$ and ending in ${\cal{D}}_j$. Obviously
that the directed graph ${\cal{D}}^*$ can not contain any cycles
since otherwise the number of strongly connected components of
${\cal{D}}$ will be equal to or less than $m-1$. It follows that
there is a strongly connected component, ${\cal{D}}_1$ say, such
that any arc that ends on a vertex in it must start at a vertex in
it. Since ${\cal{D}}$ is weakly connected, there is at least one
arc that starts in ${\cal{D}}_1$ and ends on a vertex not in
${\cal{D}}_1$. Consequently, in ${\cal{D}}_1$, the sum of
in-degree of all vertices is less than the sum of out-degree of
all vertices. This means that there must be a vertex in
${\cal{D}}$ such that the weight balance condition can not be
satisfied. Thus we have the contradiction. \hfill $\square $

Hence, if a weighted directed graph is weakly connected and the
weights of each agent satisfy the weight balance condition, then
the directed graph must be strongly connected.

\subsubsection{Stability Analysis}

Before presenting the main results of this paper, we first prove
the following important lemma.

\textbf{Lemma 1}: Let $A\in R^{n\times n}$ be any diagonal matrix
with positive diagonal entries. Then \[\left(A
{\mathrm{span}}\{\bf 1\}^\bot\right)\cap{\mathrm{span}}\{{\bf
1}\}=0,\] where ${\bf 1}=(1, \cdots, 1)^T\in{R^n}$,
${\mathrm{span}}\{{\bf 1}\}$ is the space spanned by vector ${\bf
1}$, and ${\mathrm{span}}\{\bf 1\}^\bot$ is the orthogonal
complement space of ${\mathrm{span}}\{{\bf 1}\}$.

\textbf{Proof}: Let $p\in{\left(A {\mathrm{span}}\{\bf
1\}^\bot\right)\cap{\mathrm{span}}\{\bf 1\}}.$ Then
$p\in{\mathrm{span}\{\bf 1\}}$ and there is some $q\in{
{\mathrm{span}}\{\bf 1\}^\bot}$ such that $p=Aq$. It follows that
$q^TAq=q^Tp=0.$ Since $A$ is positive definite by assumption, we
have $q=0$ and hence $p=0$. \hfill  $\square$

\textbf{Theorem 1}: By taking the control law in (\ref{eq6}),
under Assumption 2 and the weight balance condition, all agent
velocities in the group described in (\ref{eq1}) become
asymptotically the same, collision avoidance can be ensured
between all agents and the group final configuration minimizes all
agent global potentials.

\textbf{Proof}: Choose the following positive semi-definite
function
$$
J=\displaystyle\frac{1}{2}\sum_{i=1}^N\left(V^i+m_iv^{iT}v^i\right).
$$
It is easy to see that $J$ is the sum of the total artificial
potential energy and the total kinetic energy of all agents in the
group. Define the level sets of $J$ in the space of agent
velocities and relative distances
\begin{equation}
\Omega=\big\{(v^i, x^{ij}) | J\leq c \big\}.
 \label{eq7}
\end{equation}
In what follows, we will prove that the set $\Omega$ is compact.
In fact, the set $\{v^i, x^{ij}\}$ such that $J\leq c$ ($c>0$) is
closed by continuity. Moreover, boundedness can be proved under
Assumption 1, namely, from $J\leq c$, we have that $V^{ij}\leq c$.
Potential $V^{ij}$ is radially unbounded, so there must be a
positive constant $d$ such that $\|x^{ij}\|\leq d$, for all $i,
j=1, \cdots, N$. In the same way, $v^{iT}v^i\leq 2c/m_i$, thus
$\|v^i\|\leq{\sqrt{2c/m_i}}$.

By the symmetry of $V^{ij}$ with respect to $x^{ij}$ and
$x^{ij}=-x^{ji}$, it follows that
\begin{equation}
\frac{\partial V^{ij}}{\partial x^{ij}}=\frac{\partial
V^{ij}}{\partial x^i}=-\frac{\partial V^{ij}}{\partial x^j},
 \label{eq8}
\end{equation}
and therefore
$$\displaystyle\frac{d}{dt}\sum_{i=1}^N\frac{1}{2}V^i=\sum_{i=1}^N\nabla_{x^i}V^i\cdot v^i.$$

Calculating the time derivative of $J$ along the solution of
system (\ref{eq1}), we have
\begin{equation}
\begin{array}{rl}
\dot{J}= & \! \! \! -
\displaystyle\sum\limits_{i=1}^Nv^{iT}\bigg(\sum\limits_{j\sim
i}w_{ij}(v^i-v^j)\bigg)=-v^T(L\otimes I_n)v\\=& \! \! \!
\displaystyle-\frac{1}{2}v^T\left((L+L^T)\otimes I_n\right)v,
\end{array}
\label{eq9}
\end{equation}
where $v=(v^{1T}, \cdots, v^{NT})^T$ is the stack vector of all
agent velocity vectors, $L=[l_{ij}]$ with
\begin{equation}
\begin{array}{l}
{l_{ij}}=\left\{
\begin{array}{l}
-w_{ij}, \\
\sum_{k=1,k\neq i}^Nw_{ik},%
\end{array}
\begin{array}{l}
\;i\neq j, \\
\;i=j,
\end{array}
\right.  \\
\end{array}
\label{eq10}
\end{equation}
is the Laplacian matrix of the weighted velocity neighboring
graph, and $(L+L^T)\otimes I_n$ is the Kronecker product of
$L+L^T$ and $I_n$, with $I_n$ the identity matrix of order $n$.

From the definition of matrix $L$, under the weight balance
condition, it is easy to see that $L+L^T$ is symmetric and has the
properties that every row sum is equal to 0, the diagonal elements
are positive, and all the other elements are nonpositive. By
matrix theory \cite{R. Horn and C. R. Johnson}, all eigenvalues of
$L+L^T$ are nonnegative. Hence, matrix $L+L^T$ is positive
semi-definite. By the connectivity of graph ${\cal{D}}$, we know
that $L+L^T$ is irreducible and the eigenvector associated with
the single zero eigenvalue is ${\bf 1}_N$. On the other hand, it
is known that the identity matrix $I_n$ has an eigenvalue $\mu=1$
of $n$ multiplicity and $n$ linearly independent eigenvectors
\[p^1=[1, 0, \cdots, 0]^T,\
p^2=[0, 1, 0, \cdots, 0]^T,\ \cdots,\ p^n=[0, \cdots, 0, 1]^T.\]
By matrix theory \cite{R. Horn and C. R. Johnson}, the eigenvalues
of $(L+L^T)\otimes I_n$ are nonnegative, $\lambda=0$ is an
eigenvalue of multiplicity $n$ and the associated eigenvectors are
$$q^1=[p^{1T}, \cdots, p^{1T}]^T, \cdots, q^n=[p^{nT}, \cdots, p^{nT}]^T.$$
Thus $\dot{J}\leq 0$, and $\dot{J}=0$ implies that all agents have
the same velocity vector, that is, the vector $v_k=(v_{k}^1,
\cdots, v_{k}^N)$ $(k=1, \cdots, n)$, which is composed of every
corresponding $k$th component $v_{k}^1, \cdots, v_{k}^N$ of $v^1,
\cdots, v^N$, is contained in span\{{\bf 1}\}, where ${\bf 1}=(1,
\cdots, 1)^T\in R^N$. It follows that $\dot{x}^{ij}=0$, $\forall
(i, j)\in{N\times N}$.

We use LaSalle's invariance principle \cite{H. K. Khalil} to
establish convergence of system trajectories to the largest
positively invariant subset of the set defined by $E=\{v |
\dot{J}=0\}$. In $E$, the agent velocity dynamics are
$$\dot{v}^i=\frac{1}{m_i}u^i=-\frac{1}{m_i}\sum_{j=1, j\neq i}^N\nabla_{x^i}V^{ij}=-\frac{1}{m_i}\nabla_{x^i}V^i$$ and
therefore it follows that
\begin{equation}
\dot{v}=-(M\otimes I_n)\left[
\begin{array}{cc}
\nabla_{x^1}V^1 \\
\vdots\\
\nabla_{x^N}V^N%
\end{array}%
\right] = -((MB)\otimes I_n)\left[
\begin{array}{c}
\vdots \\
\nabla_{x^{ij}}V^{ij}\\
\vdots
\end{array}%
\right], \label{eq11}
\end{equation}
where $M={\mathrm {diag}} (\frac{1}{m_1}, \cdots, \frac{1}{m_N})$,
and the matrix $B$ is the incidence matrix of the position
neighboring graph. Hence
$$\dot{v}_k=-(MB)[\nabla_{x^{ij}}V^{ij}]_k,\ \  k=1, \cdots, n.$$
Thus, $\dot{v}_k\in {\mathrm {range}}(MB)$, $k=1, \cdots, n$. By
matrix theory, we have
$${\mathrm {range}}(MB)=M{\mathrm {range}}B=M{\mathrm {range}}(BB^T)=M{\mathrm {span}}\{{\bf 1}\}^\bot$$
and therefore
\begin{equation}
\dot{v}_k\in M{\mathrm {span}}\{{\bf 1}\}^\bot, \ \ k=1, \cdots,
n. \label{eq12}
\end{equation}
In any invariant set of $E$, by $v_k\in {\mathrm {span}}\{{\bf
1}\}$, we have
\begin{equation}
\dot{v}_k\in {\mathrm {span}}\{{\bf 1}\}. \label{eq13}
\end{equation}
By Lemma 1, we get from (\ref{eq12}) and (\ref{eq13})
$$\dot{v}_k\in (M{\mathrm {span}}\{{\bf 1}\}^\bot) \cap
 {\mathrm {span}}\{{\bf 1}\}\equiv {\bf 0}, \ \ k=1, \cdots, n.$$
Thus, in steady state, all agent velocities no longer change and
from (\ref{eq11}), the potential $V^i$ of each agent is globally
minimized. Collision-free can be ensured between the agents since
otherwise it will result in $V^i\rightarrow\infty$. \hfill
$\square $

\textbf{Remark 1}: If we take the control law for agent $i$ to be
\begin{equation}
u^i=-\sum_{j\in {\cal{N}}_i}(v^i-v^j)-\sum_{j=1,j\neq
i}^N\nabla_{x^i}V^{ij},\label{eq14}
\end{equation}
then the weight balance condition implies that, in the velocity
neighboring graph, for each vertex, the number of arcs starting at
it is equal to the number of arcs ending on it. When we take the
control law in (\ref{eq14}), by using the same analysis method as
in Theorem 1, we can also obtain the same conclusion.

Note that, from (\ref{eq9}), we see that the interaction
coefficients in control law (\ref{eq6}) can influence the decaying
rate of the total energy $J$. Hence, we conclude that the
convergence rate of the system will be influence by the
interaction coefficients. Explicit analysis  on this topic will be
presented in Section 4.1.3.

\subsubsection{Common Velocity}

In this section, we will show that the final common velocity can
be obtained by the initial velocities of all agents.

The position vector of the center of mass in system (\ref{eq1}) is
defined as \[x^*=\frac{\sum_{i=1}^Nm_ix^i}{\sum_{i=1}^Nm_i}.\]
Thus, the velocity vector of the center of mass is
\[v^*=\frac{\sum_{i=1}^Nm_iv^i}{\sum_{i=1}^Nm_i}.\] By using control law
(\ref{eq6}), we obtain
\begin{equation*}
\begin{array}{rl}
\dot{v}^*=
\!\!\!\displaystyle\frac{-1}{(\sum_{i=1}^Nm_i)}\sum_{i=1}^N\bigg[\sum_{j\in
{\cal{N}}_i}w_{ij}(v^i-v^j)+\sum_{j=1, j\neq i}^N\nabla_{x^i}V^{ij}\bigg].%
\end{array}
\end{equation*}%
By the symmetry of function $V^{ij}$ with respect to $x^{ij}$,
under the weight balance condition, we get $\dot{v}^*=0$. This
means that, by using control law (\ref{eq6}), the velocity of the
center of mass is invariant.

Therefore, combining Theorem 1 and the analysis above, we have the
following theorem.

\textbf{Theorem 2}: By taking the control law in (\ref{eq6}),
under Assumption 2 and the weight balance condition, the final
common velocity is equal to the initial velocity of the center of
mass, that is, the final velocity $v_{f}$
is\[v_{f}=\frac{\sum_{i=1}^Nm_iv^i(0)}{\sum_{i=1}^Nm_i},\] where
$v^i(0)$ is the velocity value of agent $i$ at initial time $t=0$,
$i=1, \cdots, N$.

\textbf{Remark 2}: Note that, by the calculation above, we can see
that the final common velocity is determined by the masses and the
initial velocities of all agents, and does not rely on the
neighboring relations and the magnitudes of the interaction
coefficients under Assumption 2 and the weight balance condition.

\textbf{Remark 3}: Even if the velocity neighboring graph is not
connected, under the weight balance condition, the velocity of the
center of mass is still invariant by using control law
(\ref{eq6}). However, in this case, the final velocities of all
agents might be different. In fact, when the velocity neighboring
graph is not connected, under the weight balance condition,
control law (\ref{eq6}) only ensures that all agents from the same
connected group will have the same final velocity, and the final
velocities of any two different connected groups might not be
equal to each other.

\textbf{Remark 4}: Using the control law in (\ref{eq6}), from
Theorems 1 and 2, we know that if the initial velocity of the
center of mass is zero, the center of mass will not drift. All
agents adjust their positions and velocities to minimize the total
potential, and the final common velocity of all agents is zero.

Hence, by using control law (\ref{eq6}), under Assumption 2, the
whole group can move ahead at a common nonzero velocity if and
only if the initial velocity of the center of mass is not zero.

\textbf{Definition 4}: The average velocity of all agents is
defined as $\overline{v}=(\sum_{i=1}^Nv^i)/N.$

\textbf{Remark 5}: If we modify the control law $u^i$ to be
\begin{equation}
u^i=-\sum_{j\in {\cal{N}}_i}m_iw_{ij}(v^i-v^j)-\sum_{j=1, j\neq
i}^Nm_i\nabla_{x^i}V^{ij},
 \label{eq15}
\end{equation}
where $m_i$ and $w_{ij}$ are defined as before, by choosing the
Lyapunov function
$$
J=\displaystyle\frac{1}{2}\sum_{i=1}^N(V^i+v^{iT}v^i),
$$
under Assumption 2 and the weight balance condition, we can still
get the results as in Theorem 1. Since the proof is similar to the
proof of Theorem 1, we omit the details.

Moreover, by using the control law in (\ref{eq15}), under
Assumption 2 and the weight balance condition, we can obtain that
the average velocity of all agents in group (\ref{eq1}) is
invariant and therefore the final velocity of the group is the
average of the initial velocities of all agents, that is,
\[v_{f}=\frac{\sum_{i=1}^Nv^i(0)}{N},\] where
$v^i(0)$ is the velocity value of agent $i$ at initial time $t=0$,
$i=1, \cdots, N$. The final common velocity does not rely on the
agents' masses, the neighboring relations, or the magnitudes of
the interaction coefficients under Assumption 2 and the weight
balance condition.

\subsubsection{Convergence Rate Analysis}

From the discussion above, we know that the coupling coefficients
can influence the decaying rate of the energy function $J$, hence,
we guess that the coupling coefficients can also influence the
convergence rate of system (\ref{eq1}). In the following, we will
present qualitative analysis of the influence of the weights
$w_{ij}$ on the convergence rate of the system.

We consider the dynamics of the error system. From the discussion
in 4.1.2, the velocity of the center of mass in system (\ref{eq1})
is invariant. Thus, we define the following error vectors:
$$e^i=x^i-x^*,$$
$$e_{v}^i=v^i-v^*,$$
where $x^*$ and $v^*$ are the position vector and the velocity
vector of the center of mass, respectively. Hence, the error
dynamics is given by
\begin{equation}
\begin{array}{l}
\dot{e}^i=e_{v}^i, \\
\dot{e}_{v}^i=\frac{1}{m_i}u^i,\ \ i=1, \cdots, N.\\
\end{array}
\label{eq16}
\end{equation}
By the definition of $V^{ij}$ and $e^i=x^i-x^*$, we get
$$\nabla_{x^i}V^{ij}(\|x^{ij}\|)=\nabla_{e^i}V^{ij}(\|e^{ij}\|).$$
By using the control law in (\ref{eq6}), we obtain
\begin{equation}
\dot{e}_{v}^i=\displaystyle\frac{1}{m_i}\bigg[-\sum\limits_{j\in{{\cal{N}}_i}}w_{ij}(e_{v}^i-e_{v}^j)-\sum\limits_{j=1, j\neq i}^N\nabla_{e^i}V^{ij}(\|e^{ij}\|)\bigg].\\
\label{eq17}
\end{equation}
We choose the following positive semi-definite function
$$J^*=\frac{1}{2}\sum_{i=1}^N\left((V^*)^i+m_ie_{v}^{iT}e_{v}^i\right)$$
which is the energy function of the error system (\ref{eq16}).
$(V^*)^i$ is the potential of agent $i$ in (\ref{eq16}) and it
equals $V^i$ by the definition of potential function $V^{ij}$.

Calculating the time derivative of $J^*$, we have
\begin{equation}
\begin{array}{rl}
\dot{J}^* = & \! \! \!
-\displaystyle\sum\limits_{i=1}^N\sum\limits_{j\in{{\cal{N}}_i}}w_{ij}e_{v}^{iT}(e_{v}^i-e_{v}^j)=-e_{v}^T(L\otimes
I_n)e_{v}\\
=& \! \! \!
-\displaystyle\frac{1}{2}e_{v}^T\left((L+L^T\right)\otimes
I_n)e_{v},
\end{array}
\label{eq18}
\end{equation}
where $e_v=(e_{v}^{1T}, \cdots, e_{v}^{NT})^T$, and $L$ and $I_n$
are defined as before.

Using the same analysis method as in Theorem 1, we have
$\dot{J}^*\leq 0$, and $\dot{J}^*=0$ implies that
$e_{v}^1=e_{v}^2=\cdots=e_{v}^N$. This occurs only when
$e_{v}^1=e_{v}^2=\cdots=e_{v}^N=0$, that is, this occurs only when
all agents have the same velocity. In other words, if there exist
two agents with different velocities, the energy function $J^*$ is
strictly monotone decreasing with time. Certainly, before the
group forms the final tight configuration, there might be the case
that all agents have the same velocity, but due to the regulation
of the potentials among agents, it instantly changes into the case
that not all agents have the same velocity except when the group
has achieved the final stable state. Hence, the decaying rate of
energy is equivalent to the convergence rate of the system. It is
easy to see that when all agents have not achieved the common
velocity, for any solution of the error system (\ref{eq16}), $e_v$
must be in the subspace spanned by eigenvectors of $(L+L^T)\otimes
I_n$ corresponding to the nonzero eigenvalues. Thus, from
(\ref{eq18}), we have $\dot{J}^*\leq -\lambda_2e_{v}^Te_v$, where
$\lambda_2$ denotes the second smallest real eigenvalue of matrix
$L+L^T$. Therefore, we have the following conclusion: The
convergence rate of the system relies on the second smallest real
eigenvalue of matrix $L+L^T$ with $L$ defined as in (\ref{eq10}).

\subsection{ Nonideal case}

We know that, in some cases, the velocity damping should not be
ignored. Then, if we still take control law (\ref{eq6}), what will
be the motion of the group? In fact, in this case, the total force
acting on the $i$th agent is
\begin{equation}
u^i=-\sum_{j\in {\cal{N}}_i}w_{ij}(v^i-v^j)-\sum_{j=1, j\neq
i}^N\nabla_{x^i}V^{ij}-k_iv^i,
 \label{eq19}
\end{equation}
where $w_{ij}$ and $k_i$ are defined as before.

The following theorem shows the motion and the final configuration
of the group.

\textbf{Theorem 3}: By taking the control law in (\ref{eq6}),
under Assumption 2 and the weight balance condition, all agent
velocities in the group described in (\ref{eq3}) become
asymptotically the same, all agents finally stop moving, collision
avoidance can be ensured between all agents, and the group final
configuration minimizes all agent global potentials.

\textbf{Proof}: Taking the Lyapunov function $J$ defined as in
Theorem 1, that is,
$$
J=\displaystyle\frac{1}{2}\sum_{i=1}^N(V^i+m_iv^{iT}v^i).
$$
We can show analogously that the set $\Omega=\{(v^i, x^{ij})|
J\leq c \}$ ($c>0$) is compact.

Calculating the time derivative of $J$, we have
$$\dot{J}=-\frac{1}{2}v^T\left((L+L^T)\otimes
 I_n\right)v-v^T\left(H\otimes
 I_n\right)v,$$
where $v$ and $L$ are defined as in Theorem 1, and
$H={\mathrm{diag}}(k_1, \cdots, k_N)$ with $k_i>0$ is the velocity
damping gain. It is easy to see that $H$ is positive definite.
Using the same analysis method as in Theorem 1, we know that
$\dot{J}\leq 0$, and $\dot{J}=0$ implies that $v^1=\cdots=v^N$ and
they all must equal zero. We denote $E^*=\{v | \dot{J}=0\}$. In
$E^*$, the agent velocity dynamics become
$$\dot{v}^i=-\frac{1}{m_i}\sum_{j=1, j\neq i}^N\nabla_{x^i}V^{ij}=-\frac{1}{m_i}\nabla_{x^i}V^i.$$
Following the proof of Theorem 1, we can conclude that
$\dot{v}_k={\bf 0}$, hence $\dot{v}^i={\bf 0}$, $i=1,\cdots,N$,
which means that the agent velocity no longer changes in steady
state. All agents will finally stop moving, and the final
configuration minimizes all agent global potentials. Furthermore,
during the course of motion, collisions can be avoided between the
agents. \hfill $\square $

\textbf{Remark 6}: It can be shown that if we use control law
(\ref{eq15}), we can still obtain all results in Theorem 3.

\textbf{Remark 7}: From Theorem 3, we know that due to damping,
all agents eventually stop moving. This is because when all agents
eventually move ahead at a common velocity, control input
(\ref{eq6}) equals zero.

In order to make the group have the same properties as in ideal
case, the control laws should contain the velocity damping term.
Hence, we modify the control scheme to be (\ref{eq4}), where
$\alpha^i$ and $\beta^i$ are defined as in (\ref{eq6}). Then, the
actual total force acting on agent $i$ is
$$
u^i=-\sum_{j\in {\cal{N}}_i}w_{ij}(v^i-v^j)-\sum_{j=1, j\neq
i}^N\nabla_{x^i}V^{ij}.
$$
Following Theorems 1 and 2, we can easily obtain the same stable
flocking motion and the final common velocity, that is, when the
velocity damping is taken into account, by using control scheme
(\ref{eq4}), under Assumption 2 and the weight balance condition,
all agent velocities in the group described in (\ref{eq3}) become
asymptotically the same, collision-free can be ensured between all
agents, the group final configuration minimizes all agent global
potentials, and the final common velocity is equal to the initial
velocity of the center of mass.

\section{Simulations}

In this section, we will present some numerical simulations for
the system described in (\ref{eq1}) in order to illustrate the
results obtained in the previous sections.

These simulations all are performed with ten agents moving on the
plane whose initial positions, velocities and the velocity
neighboring relations are selected randomly, but they satisfy: 1)
all initial positions are chosen within a ball of radius $R=15$[m]
centered at the origin, 2) all initial velocities are selected
with arbitrary directions and magnitudes in the range of (0,
10)[m/s], and 3) the velocity neighboring graph is connected. All
agents have different masses to each other and they are randomly
selected in the range of (0, 1)[kg].

Note that, because the position neighboring graph is complete, we
will not describe it. In the following figures, we only present
the velocity neighboring relations.

Figs. 1--6 show the results in one of our simulations, where the
control laws are taken in the form of (\ref{eq6}) with the
explicit potential function
$$V^{ij}=\frac{1}{2}\ln{\|x^{ij}\|^2}+\frac{5}{2\|x^{ij}\|^2},\ \ i, j=1, \cdots, 10.$$
The interaction coefficient matrix $W$ is generated randomly such
that $\sum\limits_{j=1}^{10}w_{ij}=\sum\limits_{j=1}^{10}w_{ji}$,
$w_{ii}=0$, and the nonzero $w_{ij}$ satisfy $0<w_{ij}<1$ for all
$i, j=1, \cdots, 10$. We run the simulation for 200 seconds.

In Figs. 1--4, the blue lines all represent the bidirectional
neighboring relations and the red lines with arrows represent the
unidirectional neighboring relations. Fig. 1 shows the group
initial state which includes the initial positions, velocities and
the velocity neighboring relations. Figs. 2 and 3 depict the
motion trajectories of all agents and the configurations of the
group, respectively, where the black solid arrow direction
represent the motion direction of the agents, and the dotted lines
represent the agent trajectories. In order to indicate the
influence of potential function on the group cohesion and
configuration, we present the group configuration in Fig. 2 at
time $t=60$s. It can be seen from Figs. 2 and 3 that, during the
course of motion, all agents regulate their positions to minimize
their potentials and regulate their velocities to become the same.
Fig. 4 shows the final steady state configuration and the common
velocity at $t=200$s. By numerical calculation, we can obtain that
all agents achieve the same velocity approximately at $t=128.92$s
and the final common velocity equals the initial velocity of the
center of mass. In Fig. 5, the star represents the initial
position of the center of mass, and it can be seen from it that
the velocity of the center of mass is invariant. Fig. 6 is the
velocity curves. The solid arrow indicates the tendency of
velocity variation. Fig. 6 distinctly demonstrates that all agent
velocities asymptotically approach the same.

Hence, numerical simulation also indicates that, by using the
control law in (\ref{eq6}), under the assumption of the
connectivity of the velocity neighboring graph and the weight
balance condition, stable flocking motion can be achieved.

For the case that the initial velocity of the center of mass is
zero, we also perform some simulations. Fig. 7 is one of them and
we run its associated simulation for 3000 seconds. In Fig. 7, the
star represents the position of the center of mass. In the
simulation, the center of mass is always stationary, the final
configuration no longer changes, the whole group does not drift,
and all agents finally stop moving.

\section{Conclusions}

In this paper, we have investigated the collective behavior of
multiple mobile agents moving in $n$-dimensional space with point
mass dynamics and introduced a set of control laws which enable
the group to generate stable flocking motion. We analyzed the
group properties in two different cases, respectively. When we
ignored the velocity damping, using a coordination control scheme,
we can make the group generate stable flocking motion. The control
laws are a combination of attractive/repulsive and alignment
forces and the control law acting on each agent relies on the
position information of all agents in the group and the velocity
information of all its neighboring agents. The control laws ensure
that all agent velocities become asymptotically the same,
collisions can be avoided between all agents, and the final tight
formation minimizes all agent global potentials. Moreover, we
analyzed the magnitude and direction of the final velocity and
showed that the final common velocity is equal to the initial
velocity of the center of mass of the system. When the velocity
damping is taken into account, in order to generate stable
flocking, we properly modified the control scheme such that the
velocity damping was cancelled by some terms in the control laws.
Finally, numerical simulations were worked out to further verify
our theoretical results.

\section{Appendix: Graph Theory Preliminaries}

In this section, we briefly summarize some basic concepts and
results in graph theory that have been used in this paper. More
comprehensive discussions can be found in \cite{C. Godsil and G.
Royle}.

A {\it undirected graph} $\cal{G}$ consists of a {\it vertex set}
${\cal{V}}=\{n_1, n_2, \cdots, n_m\}$ and an {\it edge set}
${\cal{E}}=\{(n_i, n_j): n_i, n_j \in {\cal{V}}\}$, where an {\it
edge} is an unordered pair of distinct vertices of $\cal{V}$. If
$n_i, n_j\in{\cal{V}}$, and $(n_i, n_j)\in{\cal{E}}$, then we say
that $n_i$ and $n_j$ are {\it adjacent} or {\it neighbors}, and
denote this by writing $n_j\sim n_i$. A graph is called {\it
complete} if every pair of vertices are adjacent. A {\it path of
length} $r$ from $n_i$ to $n_j$ in a undirected graph is a
sequence of $r+1$ distinct vertices starting with $n_i$ and ending
with $n_j$ such that consecutive vertices are adjacent. If there
is a path between any two vertices of ${\cal {G}}$, then ${\cal
{G}}$ is {\it connected}. In this paper, we always assume that the
graph is simple graph, which means that there is no self-loops and
each element of ${\cal{E}}$ is unique. An {\it oriented graph} is
a graph together with a particular orientation, where the {\it
orientation} of a graph ${\cal {G}}$ is the assignment of a
direction to each edge, so edge $(n_i, n_j)$ is an directed edge
(arc) from $n_i$ to $n_j$. The {\it incidence matrix} $B$ of an
oriented graph $\cal{G}$ is the $\{0, \pm 1\}$-matrix with rows
and columns indexed by the vertices and edges of $\cal{G}$,
respectively, such that the $ij$-entry is equal to 1 if edge $j$
is ending on vertex $n_i$, -1 if edge $j$ is beginning with vertex
$n_i$, and 0 otherwise. Define the {\it Laplacian matrix} of
$\cal{G}$ as $L{(\cal{G})}=BB^T.$ $L{(\cal{G})}$ is always
positive semi-definite. Moreover, for a connected graph,
$L{(\cal{G})}$ has a single zero eigenvalue, and the associated
right eigenvector is ${\bf 1}_m$.

A {\it directed graph} $\cal{D}$ consists of a {\it vertex set}
${\cal{V}}=\{n_1, \cdots, n_m\}$ and an {\it arc set}
${\cal{E}}=\{(n_i, n_j): n_i, n_j \in {\cal{V}}\}$, where an {\it
arc}, or {\it directed edge}, is an ordered pair of distinct
vertices of $\cal{V}$. In this paper, we always assume that
$n_i\neq n_j$, meaning that there is no self-loops, and assume
that each element of ${\cal{E}}$ is unique. Let
${\cal{D}}=({\cal{V}}, {\cal{E}}, {\cal{A}})$ be a weighted
directed graph. ${\cal{A}}=[a_{ij}]$ is the {\it weighted
adjacency matrix}, where $a_{ij}$ is the weight of arc $(n_i,
n_j)$, $a_{ij}\geq 0$ for all $i,j\in {\cal{I}}=\{1, \cdots, m\}$:
$i\neq j$ and $a_{ii}=0$ for all $i\in \cal{I}$. The set of {\it
neighbors} of vertex $n_i$ is defined as ${\cal{N}}_i=\{j\in
{\cal{I}}: a_{ij}>0\}$. The {\it in-degree} and {\it out-degree}
of vertex $n_i$ are, respectively, defined as
$${\mathrm {deg}}_{in}(n_i)=\sum\limits_{j=1}^ma_{ji}, \ \
{\mathrm {deg}}_{out}(n_i)=\sum\limits_{j=1}^ma_{ij}.$$ The {\it
weighted graph ${\cal {D}}$ Laplacian matrix} is defined as
$L({\cal{D}})=\Delta-\cal{A}$, where $\Delta$ is the {\it degree
matrix} of ${\cal {D}}$ which is a diagonal matrix and its $i$th
diagonal element is $\Delta_{ii}={\mathrm {deg}}_{out}(n_i)$. By
definition, $\lambda=0$ is an eigenvalue of the Laplacian matrix
$L({\cal{D}})$ and ${\bf 1}_m$ is its associated right
eigenvector. A {\it path} of length $r$ from $n_0$ to $n_r$ in a
directed graph is a sequence of $r+1$ distinct vertices starting
with $n_0$ and ending with $n_r$ such that $(n_{k-1}, n_k)$ is an
arc of ${\cal {D}}$ for $k=1, \cdots, r$. A {\it weak path} is a
sequence of $n_0, \cdots, n_r$ of distinct vertices such that for
$k=1, \cdots, r$, either $(n_{k-1}, n_k)$ or $(n_{k}, n_{k-1})$ is
an arc. A directed graph is {\it strongly connected} if any two
vertices can be joined by a path and is {\it weakly connected} if
any two vertices can be joined by a weak path.


\begin{thebibliography}{99}

\bibitem{N. E. Leonard and E. Fiorelli} N. E. Leonard and E. Fiorelli, ``Virtual leaders, artificial potentials and coordinated control of groups,"
in \it Proc. IEEE Conference on Decision and Control, \rm Orlando,
Florida USA, vol. 3, pp. 2968--2973, December 2001.

\bibitem{K. Warburton and J. Lazarus} K. Warburton and J. Lazarus,
``Tendency-distance models of social cohesion in animal groups,"
\it J. Theoretical Biology, \rm vol. 150, pp. 473--488, 1991.

\bibitem{I. Suzuki and M. Yamashita} I. Suzuki and M. Yamashita,
``Distributed anonymous mobile robots: Formation of geometric
patterns," \it SIAM J. Computing, \rm vol. 28, no. 4, pp.
1347--1363, 1999.

\bibitem{John H. Reif; H. Wang} J. H. Reif and H. Wang, ``Social potential fields: A distributed behavioral control for autonomous robots,"
\it Robotics and Autonomous Systems, \rm vol. 27, no. 3, pp.
171--194, May 1999.

\bibitem{F. Giulietti} F. Giulietti, L. Pollini, and M.
Innocenti, ``Autonomous formation flight," \it IEEE Control
Systems Magazine, \rm vol. 20, no. 6, pp. 34--44, December 2000.

\bibitem{E. Rimon} E. Rimon and D. E. Koditschek,
``Exact robot navigation using artificial potential functions,"
\it IEEE Transactions on Robotics and Automation, \rm vol. 8, no.
5, pp. 501--518, October 1992.

\bibitem{R. Bachmayer and N. E. Leonard} R. Bachmayer and N. E. Leonard,
``Vehicle networks for gradient descent in a sampled environment,"
in \it Proc. IEEE Conference on Decision and Control, \rm Las
Vegas, Nevada USA, vol. 1, pp. 112--117, December 2002.

\bibitem{R. Olfati-Saber 1} R. Olfati-Saber and R. M. Murray, ``Consensus problems in networks of agents with switching topology and time-delays,"
\it IEEE Transactions on Automatic Control, \rm vol. 49, no. 9,
pp. 1520--1533, September 2004.

\bibitem{C. W. Reynolds} C. W. Reynolds, ``Flocks, herds, and schools: A distributed behavioral
model," \it Computer Graphics (ACM SIGGRAPH '87 Conference
Proceedings), \rm vol. 21, no. 4, pp. 25--34, July 1987.

\bibitem{T. Vicsek} T. Vicsek, A. Czir\'{o}k, E. Ben-Jacob,
I. Cohen, and O. Shochet, ``Novel type of phase transition in a
system of self-driven particles," \it Physical Review Letters, \rm
vol. 75, no. 6, pp. 1226--1229, August 1995.

\bibitem{A. Jadbabaie} A. Jadbabaie, J. Lin, and A. S. Morse,
``Coordination of groups of mobile autonomous agents using nearest
neighbor rules," \it IEEE Transactions on Automatic Control, \rm
vol. 48, no. 6, pp. 988--1001, June 2003.

\bibitem{A. V. Savkin} A. V. Savkin,``Coordinated collective motion of
autonumous mobile robots: Analysis of Vicsek's model," \it IEEE
Transactions on Automatic Control, \rm vol. 49, no. 6, pp.
981--983, June 2004.

\bibitem{tanner 1} H. G. Tanner, A. Jadbabaie, and G. J. Pappas,
``Stable flocking of mobile agents, Part I: Fixed topology," in
\it Proc. IEEE Conference on Decision and Control, \rm Maui,
Hawaii USA, vol. 2, pp. 2010--2015, December 2003.

\bibitem{tanner 2} H. G. Tanner, A. Jadbabaie, and G. J. Pappas, ``Stable flocking of mobile agents, Part II: Dynamic topology," in
\it Proc. IEEE Conference on Decision and Control, \rm Maui,
Hawaii USA, vol. 2, pp. 2016--2021, December 2003.

\bibitem{H. Shi} H. Shi, L. Wang, T. Chu, and W. Zhang, ``Coordination of a group of mobile autonomous
agents," \it International Conference on Advances in Intelligent
Systems---Theory and Applications, \rm Luxembourg, November 2004.

\bibitem{V. Gazi and K. M. Passino} V. Gazi and K. M. Passino, ``Stability analysis of swarms,"
\it IEEE Transactions on Automatic Control, \rm vol. 48, no. 4,
pp. 692--697, April 2003.

\bibitem{V. Gazi 2} V. Gazi and K. M. Passino,
``Stability analysis of social foraging swarms," \it IEEE
Transactions on Systems, Man and Cybernetics---Part B:
Cybernetics, \rm vol. 34, no. 1, pp. 539--557, February 2004.

\bibitem{Y. Liu 2} Y. Liu, K. M. Passino, and M. Polycarpou,
``Stability analysis of m-dimensional asynchronous swarms with a
fixed communication topology," \it IEEE Transactions on Automatic
Control, \rm vol. 48, no. 1, pp. 76--95, January 2003.

\bibitem{Y. F. Liu and K. M. Passino} Y. F. Liu and K. M. Passino, ``Stable social forging swarms in a noisy environment,"
\it IEEE Transactions on Automatic Control, \rm vol. 49, no. 1,
pp. 30--44, January 2004.

\bibitem{Wang} L. Wang, H. Shi, T. Chu, W. Zhang and L. Zhang, ``Aggregation of forging swarms,"
\it Lecture Notes in Artificial Intelligence, \rm vol. 3339, pp.
766--777, Springer-Verlag, 2004.

\bibitem{T. Chu} T. Chu, L. Wang, and S. Mu, ``Collective behavior analysis of an anisotropic swarm model," in
\it Proc. of the 16th International Symposium on Mathematical
Theory of Networks and Systems, \rm Leuven, Belgium, pp. 1--14,
July 2004.

\bibitem{T. Chu2} H. Shi, L. Wang, and T. Chu, ``Swarming behavior of multi-agent systems," in
\it Proc. of the 23rd Chinese Control Conference, \rm Wuxi, China,
pp. 1027--1031, August, 2004.

\bibitem{T. Chu3} B. Liu, T. Chu, L. Wang, and F. Hao, ``Self-organization in a group of mobile autonomous agents," in
\it Proc. of the 23rd Chinese Control Conference, \rm Wuxi, China,
pp. 45--49, August, 2004.

\bibitem{T. Chu4} B. Liu, T. Chu, L. Wang., and Z. Wang, ``Swarm dynamics of a group of
mobile autonomous agents," \it Chinese Physics Letters, \rm vol.
22, no. 1, pp. 254-257, 2005.

\bibitem{R. Horn and C. R. Johnson} R. A. Horn and C. R. Johnson, \it Matrix
Analysis. \rm New York: Cambridge Univ. Press, 1985.

\bibitem{H. K. Khalil} H. K. Khalil, \it Nonlinear Systems, \rm Upper Saddle River, NJ: Prentice--Hall, 1996.

\bibitem{C. Godsil and G. Royle} C. Godsil and G. Royle, \it Algebraic Graph Theory.
\rm New York: Springer--Verlag, 2001.

\end{thebibliography}
\end{document}